\documentclass[11pt]{amsart}
\usepackage{amssymb}
\usepackage{amsmath}
\usepackage[active]{srcltx}
\usepackage{t1enc}
\usepackage[latin2]{inputenc}
\usepackage{verbatim}
\usepackage{amsmath,amsfonts,amssymb,amsthm}
\usepackage[mathcal]{eucal}
\usepackage{enumerate}
\usepackage[centertags]{amsmath}
\usepackage{graphics}

\newtheorem{theorem}{Theorem}

\newtheorem{lemma}{Lemma}

\newtheorem{corollary}{Corollary}

\begin{document}
\author{George Tephnadze}
\title[partial sums]{On the partial sums of Vilenkin-Fourier series}
\address{G. Tephnadze, Department of Mathematics, Faculty of Exact and
Natural Sciences, Tbilisi State University, Chavchavadze str. 1, Tbilisi
0128, Georgia}
\email{giorgitephnadze@gmail.com}
\thanks{The research was supported by Shota Rustaveli National Science
Foundation grant no.13/06 (Geometry of function spaces, interpolation and
embedding theorems}
\date{}
\maketitle

\begin{abstract}
The main aim of this paper is to investigate weighted maximal operators of
partial sums of Vilenkin-Fourier series. We also use our results to prove
approximation and strong convergence theorems on the martingale Hardy spaces
$H_{p},$ when $0<p\leq 1.$
\end{abstract}

\date{}



\section{ Introduction}

It is well-known that Vilenkin system forms not basis in the space $%
L_{1}\left( G_{m}\right) .$ Moreover, there is a function in the martingale
Hardy space $H_{1}\left( G_{m}\right) ,$ such that the partial sums of $f$
are not bounded in $L_{1}\left( G_{m}\right) $-norm, but partial sums $S_{n}$
of the Vilenkin-Fourier series of a function $f\in L_{1}\left( G_{m}\right) $
convergence in measure \cite{gol}.

Uniform convergence and some approximation properties of partial sums in $%
L_{1}\left( G_{m}\right) $ norms was investigate by Goginava \cite{gog1}
(see also \cite{gog2}). Fine \cite{fi} has obtained sufficient conditions
for the uniform convergence which are in complete analogy with the
Dini-Lipschits conditions. Guličev \cite{9} has estimated the rate of
uniform convergence of a Walsh-Fourier series using Lebesgue constants and
modulus of continuity. Uniform convergence of subsequence of partial sums
was investigate also in \cite{gt}. This problem has been considered for
Vilenkin group $G_{m}$ by Fridli \cite{4}, Blahota \cite{2} and Gát \cite{5}.

It is also known that subsequence $S_{n_{k}}$ is bounded from $L_{1}\left(
G_{m}\right) $ to $L_{1}\left( G_{m}\right) $ if and only if $n_{k}$ has
uniformly bounded variation and subsequence of partial sums $S_{M_{n}}$ is
bounded from the martingale Hardy space $H_{p}\left( G_{m}\right) $ to the
Lebesgue space $L_{p}\left( G_{m}\right) ,$ for all $p>0.$ In this paper we
shall prove very unexpected fact:

There exists a martingale $f\in H_{p}\left( G_{m}\right) \left( 0<p<1\right)
,$ such that
\begin{equation*}
\underset{n\in \mathbb{N}}{\sup }\left\Vert S_{M_{n}+1}f\right\Vert
_{L_{p,\infty }}=\infty .
\end{equation*}%
The reason of divergence of $S_{M_{n}+1}f$ is that when $0<p<1$ the Fourier
coefficients of $f\in H_{p}\left( G_{m}\right) $ are not bounded (See \cite%
{tep2}).

In Gát \cite{gat1} the following strong convergence result was obtained for
all  $f\in H_{1}\left( G_{m}\right) :$%
\begin{equation*}
\underset{n\rightarrow \infty }{\lim }\frac{1}{\log n}\overset{n}{\underset{%
k=1}{\sum }}\frac{\left\Vert S_{k}f-f\right\Vert _{1}}{k}=0,
\end{equation*}%
where $S_{k}f$ denotes the $k$-th partial sum of the Vilenkin-Fourier series
of $f.$ (For the trigonometric analogue see Smith \cite{sm}, for the Walsh
system see Simon \cite{Si3}). For the Vilenkin system Simon \cite{si1}
proved that there is an absolute constant $c_{p},$ depending only on $p,$
such that
\begin{equation}
\overset{\infty }{\underset{k=1}{\sum }}\frac{\left\Vert S_{k}f\right\Vert
_{p}^{p}}{k^{2-p}}\leq c_{p}\left\Vert f\right\Vert _{H_{p}}^{p},
\label{1cc}
\end{equation}%
for all $f\in H_{p}\left( G_{m}\right) $, where $0<p<1.$ The author \cite%
{tep3} proved that for any nondecreasing function $\Phi :\mathbb{N}%
\rightarrow \lbrack 1,$ $\infty )$, satisfying the condition $\underset{%
n\rightarrow \infty }{\lim }\Phi \left( n\right) =+\infty ,$ there exists a
martingale $f\in H_{p}\left( G_{m}\right) ,$ such that

\begin{equation}
\text{ }\underset{k=1}{\overset{\infty }{\sum }}\frac{\left\Vert
S_{k}f\right\Vert _{L_{p,\infty }}^{p}\Phi \left( k\right) }{k^{2-p}}=\infty
,\text{ for }0<p<1.  \label{2c}
\end{equation}

Strong convergence theorems of two-dimensional partial sums was investigate
by Weisz \cite{We}, Goginava \cite{gg}, Gogoladze \cite{Go}, Tephnadze \cite%
{tep4}.

The main aim of this paper is to investigate weighted maximal operators of
partial sums of Vilenkin-Fourier series. We also use this results to prove
some approximation and strong convergence theorems on the martingale Hardy
spaces $H_{p}\left( G_{m}\right) ,$ when $0<p\leq 1.$

\section{Definitions and Notations}

Let $\mathbb{N}_{+}$ denote the set of the positive integers, $\mathbb{N}:=%
\mathbb{N}_{+}\cup \{0\}.$

Let $m:=(m_{0,}$ $m_{1},...)$ denote a sequence of the positive integers not
less than 2.

Denote by
\begin{equation*}
Z_{m_{k}}:=\{0,1,...,m_{k}-1\}
\end{equation*}%
the additive group of integers modulo $m_{k}.$

Define the group $G_{m}$ as the complete direct product of the group $%
Z_{m_{j}}$ with the product of the discrete topologies of $Z_{m_{j}}$`s.

The direct product $\mu $ of the measures
\begin{equation*}
\mu _{k}\left( \{j\}\right) :=1/m_{k}\text{ \qquad }(j\in Z_{m_{k}})
\end{equation*}
is the Haar measure on $G_{m_{\text{ }}}$with $\mu \left( G_{m}\right) =1.$

If the sequence $m:=(m_{0,}$ $m_{1},...)$ is bounded than $G_{m}$ is called
a bounded Vilenkin group, else it is called an unbounded one.

The elements of $G_{m}$ represented by sequences
\begin{equation*}
x:=(x_{0},x_{1},...,x_{j},...)\qquad \left( \text{ }x_{k}\in
Z_{m_{k}}\right) .
\end{equation*}

It is easy to give a base for the neighborhood of $G_{m}$
\begin{equation*}
I_{0}\left( x\right) :=G_{m},
\end{equation*}%
\begin{equation*}
I_{n}(x):=\{y\in G_{m}\mid y_{0}=x_{0},...,y_{n-1}=x_{n-1}\}\text{ }(x\in
G_{m},\text{ }n\in \mathbb{N}).
\end{equation*}%
Denote $I_{n}:=I_{n}\left( 0\right) $ for $n\in \mathbb{N}$ and $\overline{%
I_{n}}:=G_{m}$ $\backslash $ $I_{n}$.

It is evident

\begin{equation}
\overline{I_{N}}=\overset{N-1}{\underset{s=0}{\bigcup }}I_{s}\backslash
I_{s+1}.  \label{2}
\end{equation}

If we define the so-called generalized number system based on $m$ in the
following way
\begin{equation*}
M_{0}:=1,\text{ \qquad }M_{k+1}:=m_{k}M_{k\text{ }}\ \qquad (k\in \mathbb{N})
\end{equation*}%
then every $n\in \mathbb{N}$ can be uniquely expressed as $n=\overset{\infty
}{\underset{k=0}{\sum }}n_{j}M_{j},$ where $n_{j}\in Z_{m_{j}}$ $~(j\in
\mathbb{N})$ and only a finite number of $n_{j}$`s differ from zero. Let $%
\left\vert n\right\vert :=\max $ $\{j\in \mathbb{N},$ $n_{j}\neq 0\}.$

Denote by $L_{1}\left( G_{m}\right) $ the usual (one dimensional) Lebesgue
space.

Next, we introduce on $G_{m}$ an orthonormal system which is called the
Vilenkin system.

At first define the complex valued function $r_{k}\left( x\right)
:G_{m}\rightarrow
\mathbb{C}
,$ the generalized Rademacher functions as
\begin{equation*}
r_{k}\left( x\right) :=\exp \left( 2\pi ix_{k}/m_{k}\right) \text{ \qquad }%
\left( \iota ^{2}=-1,\text{ }x\in G_{m},\text{ }k\in \mathbb{N}\right) .
\end{equation*}

Now define the Vilenkin system $\psi :=(\psi _{n}:n\in \mathbb{N})$ on $%
G_{m} $ as:
\begin{equation*}
\psi _{n}(x):=\overset{\infty }{\underset{k=0}{\Pi }}r_{k}^{n_{k}}\left(
x\right) \text{ \qquad }\left( n\in \mathbb{N}\right) .
\end{equation*}

Specifically, we call this system the Walsh-Paley one if $m\equiv 2$.

The Vilenkin system is orthonormal and complete in $L_{2}\left( G_{m}\right)
\,$\cite{AVD,Vi}.

Now we introduce analogues of the usual definitions in Fourier-analysis.

If $f\in L_{1}\left( G_{m}\right) $ we can establish the Fourier
coefficients, the partial sums of the Fourier series, the Dirichlet kernels
with respect to the Vilenkin system $\psi $ in the usual manner:
\begin{eqnarray*}
\widehat{f}\left( k\right) &:&=\int_{G_{m}}f\overline{\psi }_{k}d\mu ,\text{%
\thinspace }\left( \text{ }k\in \mathbb{N}\right) , \\
S_{n}f &:&=\sum_{k=0}^{n-1}\widehat{f}\left( k\right) \psi _{k},\left( \text{
}n\in \mathbb{N}_{+},\text{ }S_{0}f:=0\right) , \\
D_{n} &:&=\sum_{k=0}^{n-1}\psi _{k\text{ }},\text{ }\left( \text{ }n\in
\mathbb{N}_{+}\text{ }\right) .
\end{eqnarray*}

Recall that (see \cite{AVD})
\begin{equation}
\quad \hspace*{0in}D_{M_{n}}\left( x\right) =\left\{
\begin{array}{l}
M_{n},\text{\thinspace \thinspace \thinspace \thinspace if\thinspace
\thinspace }x\in I_{n} \\
0,\text{\thinspace \thinspace \thinspace \thinspace \thinspace if \thinspace
\thinspace }x\notin I_{n}%
\end{array}%
\right.  \label{3}
\end{equation}%
\vspace{0pt}and

\begin{equation}
D_{n}\left( x\right) =\psi _{n}(x)\left( \underset{j=0}{\overset{\infty }{%
\sum }}D_{M_{j}}\left( x\right) \overset{m_{j}-1}{\underset{u=m_{j}-n_{j}}{%
\sum }}r_{j}^{u}\left( x\right) \right) .  \label{dn5}
\end{equation}

The norm (or quasinorm) of the space $L_{p}(G_{m})$ is defined by \qquad
\qquad \thinspace\
\begin{equation*}
\left\Vert f\right\Vert _{p}:=\left( \int_{G_{m}}\left\vert f\right\vert
^{p}d\mu \right) ^{1/p}\qquad \left( 0<p<\infty \right) .
\end{equation*}%
The space $L_{p,\infty }\left( G_{m}\right) $ consists of all measurable
functions $f$ for which

\begin{equation*}
\left\Vert f\right\Vert _{L_{p},\infty }:=\underset{\lambda >0}{\sup }%
\lambda \mu \left( f>\lambda \right) ^{1/p}<+\infty .
\end{equation*}

The $\sigma $-algebra generated by the intervals $\left\{ I_{n}\left(
x\right) :x\in G_{m}\right\} $ will be denoted by $\digamma _{n}$ $\left(
n\in \mathbb{N}\right) .$ Denote by $f=\left( f_{n},n\in \mathbb{N}\right) $
a martingale with respect to $\digamma _{n}$ $\left( n\in \mathbb{N}\right)
. $ (for details see e.g. \cite{We1}). The maximal function of a martingale $%
f$ is defend by \qquad
\begin{equation*}
f^{\ast }=\sup_{n\in \mathbb{N}}\left\vert f^{\left( n\right) }\right\vert .
\end{equation*}

In case $f\in L_{1}\left( G_{m}\right) ,$ the maximal functions are also be
given by
\begin{equation*}
f^{\ast }\left( x\right) =\sup_{n\in \mathbb{N}}\frac{1}{\left\vert
I_{n}\left( x\right) \right\vert }\left\vert \int_{I_{n}\left( x\right)
}f\left( u\right) d\mu \left( u\right) \right\vert
\end{equation*}

For $0<p<\infty $ the Hardy martingale spaces $H_{p}\left( G_{m}\right) $
consist of all martingales, for which
\begin{equation*}
\left\Vert f\right\Vert _{H_{p}}:=\left\Vert f^{\ast }\right\Vert
_{p}<\infty .
\end{equation*}

The dyadic Hardy martingale spaces $H_{p}$ $\left( G_{m}\right) $ for $%
0<p\leq 1$ have an atomic characterization. Namely the following theorem is
true (see \cite{We5}):

\textbf{Theorem W}: A martingale $f=\left( f_{n},\text{ }n\in \mathbb{N}%
\right) $ is in $H_{p}\left( G_{m}\right) \left( 0<p\leq 1\right) $ if and
only if there exists a sequence $\left( a_{k},\text{ }k\in \mathbb{N}\right)
$ of p-atoms and a sequence $\left( \mu _{k},\text{ }k\in \mathbb{N}\right) $
of a real numbers such that for every $n\in \mathbb{N}$

\begin{equation}
\qquad \sum_{k=0}^{\infty }\mu _{k}S_{M_{n}}a_{k}=f_{n}  \label{2A}
\end{equation}%
and

\begin{equation*}
\qquad \sum_{k=0}^{\infty }\left\vert \mu _{k}\right\vert ^{p}<\infty .
\end{equation*}

Moreover, $\left\Vert f\right\Vert _{H_{p}}\backsim \inf \left(
\sum_{k=0}^{\infty }\left\vert \mu _{k}\right\vert ^{p}\right) ^{1/p}$,
where the infimum is taken over all decomposition of $f$ of the form (\ref%
{2A}).

Let $X=X(G_{m})$ denote either the space $L_{1}(G_{m}),$ or the space of
continuous functions $C(G_{m})$. The corresponding norm is denoted by $\Vert
.\Vert _{X}$. The modulus of continuity, when $X=C\left( G_{m}\right) $ and
the integrated modulus of continuity, where $X=L_{1}\left( G_{m}\right) $
are defined by

\begin{equation*}
\omega \left( 1/M_{n},f\right) _{X}=\sup\limits_{h\in I_{n}}\left\Vert
f\left( \cdot +h\right) -f\left( \cdot \right) \right\Vert _{X}.
\end{equation*}

The concept of modulus of continuity in $H_{p}\left( G_{m}\right) \left(
0<p\leq 1\right) $ can be defined in following way
\begin{equation*}
\omega \left( 1/M_{n},f\right) _{H_{p}\left( G_{m}\right) }:=\left\Vert
f-S_{M_{n}}f\right\Vert _{H_{p}\left( G_{m}\right) }.
\end{equation*}

If $f\in L_{1}\left( G_{m}\right) ,$ then it is easy to show that the
sequence $\left( S_{M_{n}}\left( f\right) :n\in \mathbb{N}\right) $ is a
martingale.

If $f=\left( f_{n},n\in \mathbb{N}\right) $ is martingale then the
Vilenkin-Fourier coefficients must be defined in a slightly different
manner: $\qquad \qquad $
\begin{equation*}
\widehat{f}\left( i\right) :=\lim_{k\rightarrow \infty
}\int_{G_{m}}f_{k}\left( x\right) \overline{\Psi }_{i}\left( x\right) d\mu
\left( x\right) .
\end{equation*}%
\qquad \qquad \qquad

The Vilenkin-Fourier coefficients of $f\in L_{1}\left( G_{m}\right) $ are
the same as the martingale $\left( S_{M_{n}}\left( f\right) :n\in \mathbb{N}%
\right) $ obtained from $f$ .

For the martingale $f$ we consider maximal operators
\begin{eqnarray*}
S^{\ast }f &:&=\sup_{n\in \mathbb{N}}\left\vert S_{n}f\right\vert , \\
\widetilde{S}_{p}^{\ast }f &:&=\sup_{n\in \mathbb{N}}\frac{\left\vert
S_{n}f\right\vert }{\left( n+1\right) ^{1/p-1}\log ^{\left[ p\right] }\left(
n+1\right) },\text{ }0<p\leq 1,
\end{eqnarray*}%
where $\left[ p\right] $ denotes integer part of $p.$

A bounded measurable function $a$ is p-atom, if there exist a dyadic
interval $I$, such that%
\begin{equation*}
\int_{I}ad\mu =0,\text{ \ \ \ \ }\left\Vert a\right\Vert _{\infty }\leq \mu
\left( I\right) ^{-1/p},\text{ \ \ \ \ \ supp}\left( a\right) \subset I.
\end{equation*}
\qquad

\section{Formulation of Main Results}

\begin{theorem}
a) Let $0<p\leq 1.$ Then the \bigskip maximal operator $\overset{\sim }{S}%
_{p}^{\ast }$ \textit{is bounded from the Hardy space }$H_{p}\left(
G_{m}\right) $\textit{\ to the space }$L_{p}\left( G_{m}\right) .$
\end{theorem}

b) Let $0<p\leq 1$ and $\varphi :\mathbb{N}_{+}\rightarrow \lbrack 1,$ $%
\infty )$ be a nondecreasing function satisfying the condition

\begin{equation}
\overline{\lim_{n\rightarrow \infty }}\frac{\left( n+1\right) ^{1/p-1}\log ^{%
\left[ p\right] }\left( n+1\right) }{\varphi \left( n\right) }=+\infty .
\label{6}
\end{equation}%
\textit{Then }

\begin{equation*}
\sup_{n\in \mathbb{N}}\left\Vert \frac{S_{n}f}{\varphi \left( n\right) }%
\right\Vert _{L_{p,\infty }\left( G_{m}\right) }=\infty ,\text{ for }0<p<1
\end{equation*}%
and%
\begin{equation*}
\sup_{n\in \mathbb{N}}\left\Vert \frac{S_{n}f}{\varphi \left( n\right) }%
\right\Vert _{1}=\infty .
\end{equation*}

\begin{corollary}
(Simon \cite{si1}) Let $0<p<1\ $and $f\in H_{p}\left( G_{m}\right) .$ Then
there is an absolute constant $c_{p},$ depends only $p,$ such that
\begin{equation*}
\overset{\infty }{\underset{k=1}{\sum }}\frac{\left\Vert S_{k}f\right\Vert
_{p}^{p}}{k^{2-p}}\leq c_{p}\left\Vert f\right\Vert _{H_{p}}^{p}.
\end{equation*}
\end{corollary}

\begin{theorem}
Let $0<p\leq 1,$ $f\in H_{p}\left( G_{m}\right) $ and $M_{k}<n\leq M_{k+1}$.
Then there is an absolute constant $c_{p},$ depends only $p,$ such that
\begin{equation*}
\left\Vert S_{n}\left( f\right) -f\right\Vert _{H_{p}\left( G_{m}\right)
}\leq c_{p}n^{1/p-1}\lg ^{\left[ p\right] }n\omega \left( \frac{1}{M_{k}}%
,f\right) _{H_{p}\left( G_{m}\right) }.
\end{equation*}
\end{theorem}

\begin{theorem}
a) Let $0<p<1,$ $f\in H_{p}\left( G_{m}\right) $ and%
\begin{equation*}
\omega \left( \frac{1}{M_{n}},f\right) _{H_{p}\left( G_{m}\right) }=o\left(
\frac{1}{M_{n}^{1/p-1}}\right) ,\text{ as }n\rightarrow \infty .
\end{equation*}

Then%
\begin{equation*}
\left\Vert S_{k}\left( f\right) -f\right\Vert _{L_{p,\infty }\left(
G_{m}\right) }\rightarrow 0,\,\,\,\text{when\thinspace \thinspace \thinspace
}k\rightarrow \infty .
\end{equation*}
\end{theorem}

b) For every $p\in \left( 0,1\right) $ there exists martingale $f\in
H_{p}(G_{m})$, for which
\begin{equation*}
\omega \left( \frac{1}{M_{2n}},f\right) _{H_{p}(G_{m})}=O\left( \frac{1}{%
M_{2n}^{1/p-1}}\right) ,\text{ as }n\rightarrow \infty
\end{equation*}%
and
\begin{equation*}
\left\Vert S_{k}\left( f\right) -f\right\Vert _{L_{p,\infty
}(G_{m})}\nrightarrow 0,\,\,\,\text{when\thinspace \thinspace \thinspace }%
k\rightarrow \infty
\end{equation*}

\begin{theorem}
Let $f\in H_{1}(G_{m})\ $and
\begin{equation*}
\omega \left( \frac{1}{M_{n}},f\right) _{H_{1}(G_{m})}=o\left( \frac{1}{n}%
\right) ,\text{ as }n\rightarrow \infty .
\end{equation*}%
Then%
\begin{equation*}
\left\Vert S_{k}\left( f\right) -f\right\Vert _{1}\rightarrow 0,\,\,\,\text{%
when\thinspace \thinspace \thinspace }k\rightarrow \infty .
\end{equation*}

b) There exists martingale $f\in H_{1}(G_{m})$ \ for which
\begin{equation*}
\omega \left( \frac{1}{M_{2M_{n}}},f\right) _{H_{1}(G_{m})}=O\left( \frac{1}{%
M_{n}}\right) ,\text{ as }n\rightarrow \infty
\end{equation*}%
and
\begin{equation*}
\left\Vert S_{k}\left( f\right) -f\right\Vert _{1}\nrightarrow 0\,\,\,\text{%
when\thinspace \thinspace \thinspace }k\rightarrow \infty .
\end{equation*}
\end{theorem}

\section{Auxiliary propositions}

\begin{lemma}
\cite{We3} Suppose that an operator $T$ is sublinear and for some $0<p\leq 1$
\end{lemma}

\begin{equation*}
\int\limits_{\overline{I}}\left\vert Ta\right\vert ^{p}d\mu \leq
c_{p}<\infty ,
\end{equation*}%
for every $p$-atom $a$, where $I$ denote the support of the atom. If $T$ is
bounded from $L_{\infty \text{ }}$ to $L_{\infty }.$ Then
\begin{equation*}
\left\Vert Tf\right\Vert _{p}\leq c_{p}\left\Vert f\right\Vert
_{H_{p}(G_{m})}.
\end{equation*}

\begin{lemma}
\cite{tep2} Let $n\in \mathbb{N}$ and $x\in I_{s}\backslash I_{s+1},$ $0\leq
s\leq N-1.$ Then
\end{lemma}

\begin{equation*}
\int_{I_{N}}\left\vert D_{n}\left( x-t\right) \right\vert d\mu \left(
t\right) \leq \,\frac{cM_{s}}{M_{N}}.
\end{equation*}

\section{Proof of the Theorems}

\textbf{Proof of Theorem 1. }Since $\overset{\sim }{S}_{p}^{\ast }$ is
bounded from $L_{\infty }(G_{m})$ to $L_{\infty }(G_{m})$ by Lemma 1 we
obtain that the proof of theorem 1 will be complete, if we show that

\begin{equation*}
\int\limits_{\overline{I}_{N}}\left\vert \overset{\sim }{S}_{p}^{\ast
}a\left( x\right) \right\vert ^{p}d\mu \left( x\right) \leq c<\infty ,\text{
when }0<p\leq 1,
\end{equation*}%
for every $p$-atom $a,$ where $I$ denotes the support of the atom$.$

Let $a$ be an arbitrary $p$-atom with support$\ I$ and $\mu \left( I\right)
=M_{N}.$ We may assume that $I=I_{N}.$ It is easy to see that $S_{n}\left(
a\right) =0$ when $n\leq M_{N}.$ Therefore we can suppose that $n>M_{N}.$

Since $\left\Vert a\right\Vert _{\infty }\leq M_{N}^{1/p}$ we can write
\begin{eqnarray}
&&\left\vert S_{n}\left( a\right) \right\vert \leq \int_{I_{N}}\left\vert
a\left( t\right) \right\vert \left\vert D_{n}\left( x-t\right) \right\vert
d\mu \left( t\right)  \label{12A} \\
&\leq &\left\Vert a\right\Vert _{\infty }\int_{I_{N}}\left\vert D_{n}\left(
x-t\right) \right\vert d\mu \left( t\right) \leq
M_{N}^{1/p}\int_{I_{N}}\left\vert D_{n}\left( x-t\right) \right\vert d\mu
\left( t\right) .  \notag
\end{eqnarray}

Let $0<p<1$ and $x\in I_{s}\backslash I_{s+1}.$ From Lemma 2 we get
\begin{equation}
\frac{\left\vert S_{n}a\left( x\right) \right\vert }{\log ^{\left[ p\right]
}\left( n+1\right) \left( n+1\right) ^{1/p-1}}\leq \frac{cM_{N}^{1/p-1}M_{s}%
}{\log ^{\left[ p\right] }\left( n+1\right) \left( n+1\right) ^{1/p-1}}.
\label{13AA}
\end{equation}%
Combining (\ref{2}) and (\ref{13AA}) we obtain
\begin{eqnarray}
&&\int_{\overline{I_{N}}}\left\vert \widetilde{S}_{p}^{\ast }a\left(
x\right) \right\vert ^{p}d\mu \left( x\right) =\overset{N-1}{\underset{s=0}{%
\sum }}\int_{I_{s}\backslash I_{s+1}}\left\vert \widetilde{S}_{p}^{\ast
}a\left( x\right) \right\vert ^{p}d\mu \left( x\right)  \label{14A} \\
&\leq &\frac{cM_{N}^{1-p}}{\log ^{\left[ p\right] p}\left( n+1\right) \left(
n+1\right) ^{1-p}}\overset{N-1}{\underset{s=0}{\sum }}\frac{M_{s}^{p}}{M_{s}}%
\leq \frac{cM_{N}^{1-p}N^{\left[ p\right] }}{\log ^{p\left[ p\right] }\left(
n+1\right) \left( n+1\right) ^{1-p}}<c_{p}<\infty .  \notag
\end{eqnarray}

Let $0<p<1.$ Applying (\ref{12A}), (\ref{14A}) and Theorem W we have%
\begin{eqnarray}
&&\overset{\infty }{\underset{k=M_{N}}{\sum }}\frac{\left\Vert
S_{k}a\right\Vert _{p}^{p}}{k^{2-p}}\leq \overset{\infty }{\underset{k=M_{N}}%
{\sum }}\frac{1}{k}\int_{\overline{I_{N}}}\left\vert \frac{S_{k}a\left(
x\right) }{k^{1/p-1}}\right\vert ^{p}d\mu \left( x\right)  \label{11s} \\
&&+\overset{\infty }{\underset{k=M_{N}}{\sum }}\frac{M_{N}}{k^{2-p}}%
\int_{I_{N}}\left( \int_{I_{N}}\left\vert D_{k}\left( x-t\right) \right\vert
d\mu \left( t\right) \right) ^{p}d\mu \left( x\right)  \notag \\
&\leq &c_{p}M_{N}^{1-p}\overset{\infty }{\underset{k=M_{N}}{\sum }}\frac{1}{%
k^{2-p}}+c_{p}M_{N}^{1-p}\overset{\infty }{\underset{k=M_{N}}{\sum }}\frac{%
\log ^{p}k}{k^{2-p}}\leq c_{p}<\infty .  \notag
\end{eqnarray}

Which complete the proof of corollary 1.

Let prove second part of Theorem 1. Let

\begin{equation*}
f_{n_{k}}\left( x\right) =D_{M_{2n_{k}+1}}\left( x\right)
-D_{M_{_{2n_{k}}}}\left( x\right) .\text{ }
\end{equation*}

It is evident
\begin{equation*}
\widehat{f}_{n_{k}}\left( i\right) =\left\{
\begin{array}{l}
\text{ }1,\text{ if }i=M_{_{2n_{k}}},...,M_{2n_{k}+1}-1 \\
\text{ }0,\text{otherwise.}%
\end{array}%
\right.
\end{equation*}%
Then we can write
\begin{equation}
S_{i}f_{n_{k}}\left( x\right) =\left\{
\begin{array}{l}
D_{i}\left( x\right) -D_{M_{_{2n_{k}}}}\left( x\right) ,\text{ if }%
i=M_{_{2n_{k}}}+1,...,M_{2n_{k}+1}-1, \\
\text{ }f_{n_{k}}\left( x\right) ,\text{ if }i\geq M_{2n_{k}+1}, \\
0,\text{ otherwise.}%
\end{array}%
\right.  \label{14}
\end{equation}

From (\ref{3}) we get
\begin{equation}
\left\Vert f_{n_{k}}\right\Vert _{H_{p}(G_{m})}=\left\Vert \sup\limits_{n\in
\mathbb{N}}S_{M_{n}}\left( f_{n_{k}}\right) \right\Vert _{p}=\left\Vert
D_{M_{2n_{k}+1}}-D_{M_{_{2n_{k}}}}\right\Vert _{p}\leq
c_{p}M_{_{2n_{k}}}^{1-1/p}.  \label{1}
\end{equation}%
Let $0<p<1.$ Under condition (\ref{6}) there exists positive integers $n_{k}$
such that
\begin{equation*}
\lim_{k\rightarrow \infty }\frac{\left( M_{2n_{k}}+2\right) ^{1/p-1}}{%
\varphi \left( M_{2n_{k}}+2\right) }=\infty ,\text{ \ \ }0<p<1.
\end{equation*}%
Applying (\ref{3}), (\ref{dn5}) and (\ref{14}) we can write
\begin{equation*}
\frac{\left\vert S_{M_{2n_{k}}+1}f_{n_{k}}\right\vert }{\varphi \left(
M_{2n_{k}}+2\right) }=\frac{\left\vert
D_{M_{2n_{k}}+1}-D_{M_{2_{n_{k}}}}\right\vert }{\varphi \left(
M_{2n_{k}}+2\right) }=\frac{\left\vert w_{M_{2_{n_{k}}}}\right\vert }{%
\varphi \left( M_{2n_{k}}+2\right) }=\frac{1}{\varphi \left(
M_{2n_{k}}+2\right) }.
\end{equation*}

Hence we can write:

\begin{equation}
\mu \left\{ x\in G_{m}:\frac{\left\vert S_{M_{2n_{k}}+1}f_{n_{k}}\left(
x\right) \right\vert }{\varphi \left( M_{2n_{k}}+2\right) }\geq \frac{1}{%
\varphi \left( M_{2n_{k}}+2\right) }\right\} =1.  \label{16}
\end{equation}

Combining (\ref{1}) and (\ref{16}) we have

\begin{eqnarray*}
&&\frac{\frac{1}{\varphi \left( M_{2n_{k}}+2\right) }\left( \mu \left\{ x\in
G_{m}:\frac{\left\vert S_{M_{2n_{k}}+1}f_{n_{k}}\left( x\right) \right\vert
}{\varphi \left( M_{2n_{k}}+2\right) }\geq \frac{1}{\varphi \left(
M_{2n_{k}}+2\right) }\right\} \right) ^{1/p}}{\left\Vert f_{n_{k}}\left(
x\right) \right\Vert _{H_{p}}} \\
&\geq &\frac{1}{\varphi \left( M_{2n_{k}}+2\right) M_{_{2n_{k}}}^{1-1/p}}=%
\frac{\left( M_{_{2n_{k}}}+2\right) ^{1/p-1}}{\varphi \left(
M_{2n_{k}}+2\right) }\rightarrow \infty ,\text{ when \ }k\rightarrow \infty .
\end{eqnarray*}

Now consider the case when $p=1.$ Under condition (\ref{6}) there exists $%
\left\{ n_{k}:k\geq 1\right\} ,$ such that
\begin{equation*}
\lim_{k\rightarrow \infty }\frac{\log q_{n_{k}}}{\varphi \left(
q_{n_{k}}\right) }=\infty .
\end{equation*}

Let $q_{n_{k}}$\bigskip $=M_{2n_{k}}+M_{2n_{k}-2}+M_{2}+M_{0}$ and $x\in
I_{2s}\backslash I_{2s+1},$ $s=0,...,n_{k}.$ Combining (\ref{3}) and (\ref%
{dn5}) we have
\begin{eqnarray*}
&&\left\vert D_{q_{n_{k}}}\left( x\right) \right\vert \geq \left\vert
D_{M_{2s}}\left( x\right) \right\vert -\left\vert \underset{l=0}{\overset{s-2%
}{\sum }}r_{2l}^{m_{2l}-1}\left( x\right) D_{M_{2l}}\left( x\right)
\right\vert \\
&\geq &M_{2s}-\underset{l=0}{\overset{s-2}{\sum }}M_{2l}\geq
M_{2s}-M_{2s-1}\geq \frac{M_{2s}}{2}.
\end{eqnarray*}

Hence%
\begin{equation}
\int_{G_{m}}\left\vert D_{q_{n_{k}}}\left( x\right) \right\vert d\mu \left(
x\right) \geq \frac{1}{2}\underset{s=0}{\overset{n_{k}}{\sum }}%
\int_{I_{2s}\backslash I_{2s+1}}M_{2s}d\mu \left( x\right) \geq c\underset{%
s=0}{\overset{n_{k}}{\sum }}1\geq cn_{k}.  \label{22}
\end{equation}

From (\ref{14}), (\ref{1}) and (\ref{22}) we have

\begin{eqnarray*}
&&\frac{1}{\left\Vert f_{n_{k}}\left( x\right) \right\Vert _{H_{1}(G_{m})}}%
\int_{G_{m}}\frac{\left\vert S_{q_{n_{k}}}f_{n_{k}}\left( x\right)
\right\vert }{\varphi \left( q_{n_{k}}\right) }d\mu \left( x\right) \\
&\geq &\frac{1}{\left\Vert f_{n_{k}}\left( x\right) \right\Vert
_{H_{1}(G_{m})}}\left( \int_{G_{m}}\frac{\left\vert D_{q_{n_{k}}}\left(
x\right) \right\vert }{\varphi \left( q_{n_{k}}\right) }d\mu \left( x\right)
-\int_{G_{m}}\frac{\left\vert D_{M_{_{2n_{k}}}}\left( x\right) \right\vert }{%
\varphi \left( q_{n_{k}}\right) }d\mu \left( x\right) \right) \\
&\geq &\frac{c}{\varphi \left( q_{n_{k}}\right) }\left( \log
q_{n_{k}}-1\right) \geq \frac{c\log q_{n_{k}}}{\varphi \left(
q_{n_{k}}\right) }\rightarrow \infty ,\text{ \qquad when }k\rightarrow
\infty .
\end{eqnarray*}

Which complete the proof of theorem 1.

\textbf{Proof of Theorem 2. }Let $0<p\leq 1$ and $M_{k}<n\leq M_{k+1}.$
Using Theorem 1 we have
\begin{equation*}
\left\Vert S_{n}f\right\Vert _{p}\leq c_{p}n^{1/p-1}\log ^{\left[ p\right]
}n\left\Vert f\right\Vert _{H_{p}(G_{m})}.
\end{equation*}

Hence

\begin{eqnarray*}
&&\left\Vert S_{n}f-f\right\Vert _{p}^{p}\leq \left\Vert
S_{n}f-S_{M_{k}}f\right\Vert _{p}^{p}+\left\Vert S_{M_{k}}f-f\right\Vert
_{p}^{p}=\left\Vert S_{n}\left( S_{M_{k}}f-f\right) \right\Vert _{p}^{p} \\
&&+\left\Vert S_{M_{k}}f-f\right\Vert _{p}^{p}\leq c_{p}\left(
n^{1-p}+1\right) \log ^{p\left[ p\right] }n\omega ^{p}\left( \frac{1}{M_{k}}%
,f\right) _{H_{p}(G_{m})}
\end{eqnarray*}%
and%
\begin{equation}
\left\Vert S_{n}f-f\right\Vert _{p}\leq c_{p}n^{1/p-1}\log ^{\left[ p\right]
}n\omega \left( \frac{1}{M_{k}},f\right) _{H_{p}(G_{m})}.  \label{app}
\end{equation}

\textbf{Proof of Theorem 3. }Let $0<p<1,$ $f\in H_{p}(G_{m})$ and%
\begin{equation*}
\omega \left( \frac{1}{M_{2n}},f\right) _{H_{p}(G_{m})}=o\left( \frac{1}{%
M_{2n}^{1/p-1}}\right) ,\text{ as }n\rightarrow \infty .
\end{equation*}%
Using (\ref{app}) we immediately get%
\begin{equation*}
\left\Vert S_{n}f-f\right\Vert _{p}\rightarrow \infty ,\text{ when }%
n\rightarrow \infty .
\end{equation*}

Let proof of second part of theorem 3. We set
\begin{equation*}
a_{k}\left( x\right) =\frac{M_{2k}^{1/p-1}}{\lambda }\left(
D_{M_{2k+1}}\left( x\right) -D_{M_{2k}}\left( x\right) \right) ,
\end{equation*}%
where $\lambda =\sup_{n\in \mathbb{N}}m_{n}$ and
\begin{equation*}
f_{A}\left( x\right) =\underset{i=0}{\overset{A}{\sum }}\frac{\lambda }{%
M_{2i}^{1/p-1}}a_{i}(x).
\end{equation*}

Since

\begin{equation}
S_{M_{A}}a_{k}\left( x\right) =\left\{
\begin{array}{l}
a_{k}\left( x\right) ,\text{ \ }2k\leq A, \\
0,\text{ \ \ \ \ \ \ \ \ }2k>A,%
\end{array}%
\right.  \label{19}
\end{equation}%
and

\begin{equation*}
\text{supp}(a_{k})=I_{2k},\text{ \ \ }\int_{I_{2k}}a_{k}d\mu =0,\text{ \ \ }%
\left\Vert a_{k}\right\Vert _{\infty }\leq M_{2k}^{1/p-1}\cdot
M_{2k}=M_{2k}^{1/p}=\left( \text{supp }a_{k}\right) ^{-1/p},
\end{equation*}

if we apply Theorem W we conclude that $f\in H_{p}.$

It is easy to show that
\begin{eqnarray}
&&f-S_{M_{n}}f  \label{20} \\
&=&\left( f^{\left( 1\right) }-S_{M_{n}}f^{\left( 1\right) },...,f^{\left(
n\right) }-S_{M_{n}}f^{\left( n\right) },...,f^{\left( n+k\right)
}-S_{M_{n}}f^{\left( n+k\right) }\right)  \notag \\
&=&\left( 0,...,0,f^{\left( n+1\right) }-f^{\left( n\right) },...,f^{\left(
n+k\right) }-f^{\left( n\right) },...\right)  \notag \\
&=&\left( 0,...,0,\underset{i=n}{\overset{k}{\sum }}\frac{a_{i}(x)}{%
M_{i}^{1/p-1}},...\right) ,\text{ \ }k\in \mathbb{N}_{+}  \notag
\end{eqnarray}%
is martingale. Using (\ref{20}) we get%
\begin{equation*}
\omega (\frac{1}{M_{n}},f)_{H_{p}}\leq \sum\limits_{i=\left[ n/2\right]
+1}^{\infty }\frac{1}{M_{2i}^{1/p-1}}=O\left( \frac{1}{M_{n}^{1/p-1}}\right)
.
\end{equation*}%
where $\left[ n/2\right] $ denotes integer part of $n/2.$ It is easy to show
that

\begin{equation}
\widehat{f}(j)=\left\{
\begin{array}{l}
1,\,\ \text{if \thinspace \thinspace }j\in \left\{
M_{2i},...,M_{2i+1}-1\right\} ,\text{ }i=0,1,... \\
0,\text{ \thinspace \thinspace \thinspace if \thinspace \thinspace
\thinspace }j\notin \bigcup\limits_{i=0}^{\infty }\left\{
M_{2i},...,M_{2i+1}-1\right\} .\text{ }%
\end{array}%
\right.  \label{29}
\end{equation}

Using (\ref{29}) we have
\begin{equation*}
\limsup\limits_{k\rightarrow \infty }\Vert f-S_{M_{2k+1}-1}(f)\Vert
_{L_{p,\infty }(G_{m})}
\end{equation*}%
\begin{equation*}
\geq \limsup\limits_{k\rightarrow \infty }\left( \Vert w_{M_{2k+1}-1}\Vert
_{L_{p,_{\infty }}(G_{m})}-\Vert \sum\limits_{i=k+1}^{\infty }\left(
D_{M_{2i+1}}-D_{M_{2i}}\right) \Vert _{L_{p,\infty }(G_{m})}\right)
\end{equation*}%
\begin{equation*}
\geq \limsup\limits_{k\rightarrow \infty }\left( 1-c/M_{2k}^{1/p-1}\right)
c>0.
\end{equation*}

Which complete the proof of Theorem 3.

\textbf{Proof of Theorem 4. }Analogously we can prove first part of Theorem
4.\textbf{\ }Let proof it`s second part. We set
\begin{equation*}
a_{i}(x)=D_{M_{2M_{i}+1}}(x)-D_{M_{2M_{i}}}(x)
\end{equation*}%
and
\begin{equation*}
f_{A}\left( x\right) =\overset{A}{\sum_{i=1}}\frac{a_{i}(x)}{M_{i}}.
\end{equation*}

Since

\begin{equation}
S_{M_{A}}a_{k}\left( x\right) =\left\{
\begin{array}{l}
a_{k}\left( x\right) ,\text{ }2M_{k}\leq A, \\
0,\text{ \qquad }2M_{k}>A,%
\end{array}%
\right.  \label{10}
\end{equation}%
and

\begin{equation*}
\text{supp}(a_{k})=I_{2M_{k}},\text{ \ }\int_{I_{2M_{k}}}a_{k}d\mu =0,\text{
\ }\left\Vert a_{k}\right\Vert _{\infty }\leq M_{2M_{k}}=\mu (\text{supp }%
a_{k}),
\end{equation*}

if we apply Theorem W we conclude that $f\in H_{1}.$

It is easy to show that%
\begin{equation}
\omega (\frac{1}{M_{n}},f)_{H_{1}(G_{m})}\leq \sum\limits_{i=\left[ \lg n/2%
\right] }^{\infty }\frac{1}{M_{i}}=O\left( \frac{1}{n}\right) ,  \label{12}
\end{equation}%
where $\left[ \lg n/2\right] $ denotes integer part of $\lg n/2.$ By simple
calculation we get

\begin{equation}
\widehat{f}(j)=\left\{
\begin{array}{l}
\frac{1}{M_{2i}},\,\ \text{if \thinspace \thinspace }j\in \left\{
M_{2M_{i}},...,M_{2M_{i}+1}-1\right\} ,\text{ }i=0,1,... \\
0,\text{ \thinspace \thinspace \thinspace if \thinspace \thinspace
\thinspace }j\notin \bigcup\limits_{i=0}^{\infty }\left\{
M_{2M_{i}},...,M_{2M_{i}+1}-1\right\} .\text{ }%
\end{array}%
\right.  \label{13}
\end{equation}

Combining (\ref{22}) and (\ref{13}) we have
\begin{eqnarray*}
&&\limsup\limits_{k\rightarrow \infty }\Vert f-S_{q_{M_{k}}}(f)\Vert _{1} \\
&\geq &\limsup\limits_{k\rightarrow \infty }\left( \frac{1}{M_{2k}}\Vert
D_{q_{M_{k}}}\Vert _{1}-\frac{1}{M_{2k}}\Vert D_{M_{2M_{k}+1}}\Vert
_{1}-\Vert \overset{\infty }{\sum_{i=k+1}}\frac{%
D_{M_{2M_{i}+1}}-D_{M_{2M_{i}}}}{M_{2i}}\Vert _{1}\right) \\
&\geq &\limsup\limits_{k\rightarrow \infty }\left( c-\overset{\infty }{%
\sum_{i=k+1}}\frac{1}{M_{2i}}-\frac{1}{M_{2k}}\right) \geq c>0.
\end{eqnarray*}

Theorem 4 is proved.

\end{document}